\documentclass[journal]{IEEEtran}
\ifCLASSINFOpdf
\else
\fi
\usepackage{amsmath}

\usepackage{stfloats}
\usepackage{threeparttable}
\usepackage{multirow}
\usepackage{diagbox}
\usepackage{url}
\usepackage[pdftex]{graphicx}
\graphicspath{{./figure/}}
\DeclareGraphicsExtensions{.pdf,.jpeg,.png,.jpg}

\begin{document}
%
\title{DC-based Security Constraints Formulation: A Perspective of Primal-Dual Interior Point Method}
%
%
%

\author{Zhiyuan~Bao,~\IEEEmembership{Student Member,~IEEE},
        Zechun~Hu,~\IEEEmembership{Senior Member,~IEEE}, Asad Mujeeb%
        \thanks{Z. Bao, Z. Hu, and A. Mujeeb are with the Department of Electrical Engineering, Tsinghua University, Beijing 100084, China (baozy19@mails.tsinghua.edu.cn; zechhu@tsinghua.edu.cn; asd20@mails.tsinghua.edu.cn).}
}

\maketitle

\begin{abstract}
The DC network security constraints have been extensively studied in numerous power system problems, such as optimal power flow (OPF), security-constrained economic dispatch (SCED), and security-constrained unit commitment (SCUC). Linear shift factors, i.e., power transfer distribution factors (PTDFs), are widely applied to replace DC power flow constraints. However, the PTDF matrix is extremely dense, making it difficult to solve security-constraint optimization problems. This paper analyzes/investigates the computational inefficiency of PTDF-based security constraints from the sparse structure perspective of the primal-dual interior point method(IPM). Additionally, a matrix transformation method is proposed for restoring the sparsity of the linear system during IPM iterations. It turns out that the transformation method is equivalent to solving the original optimization problem expressed in pure voltage angle, which preserves the sparsity structure but introduces additional variables and constraints proportional to one to two times the total number of buses. The regular B-$\theta$ formulation is also a variant of the proposed transformation. Numerical studies show that sparsity rather than the size of variables and constraints is the key factor impacting the speed of solving convex quadratic problems (QP), i.e., OPF and SCED problems. In contrast, sparsity is less desirable when solving a mixed integer problem (MIP), such as the SCUC problem, where reoptimization techniques are significantly more critical and the dual simplex method is typically employed rather than IPM.
\end{abstract}

\begin{IEEEkeywords}
Security constraints, primal-dual interior point method,  power transfer distribution factors, sparsity structure, economic dispatch.
\end{IEEEkeywords}

%
\IEEEpeerreviewmaketitle

\section{Introduction}
%
%
%
%
\IEEEPARstart{P}{ower} flow security constraints are the fundamental building block of numerous power system optimization problems, including Optimal Power Flow (OPF), Security-Constrained Economic Dispatch (SCED), and Security-Constrained Unit Commitment (SCUC). In most electricity markets regulated by Independent System Operators (ISOs), it is necessary to solve a real-time Security-Constrained Economic Dispatch (RT-SCED) problem within a timeframe of fewer than five minutes. Due to the gradual increase in the proportion of renewable energy and the increasing uncertainty in operation, operators must monitor the risks associated with various operational scenarios in real-time \cite{chen2022learning}. Given the narrow time window for real-time operation decision, it is necessary to conduct a fine-grained analysis of security constraints to improve solving speed.

Due to a large number of power flow security constraints, it becomes the bottleneck in the solution of SCUC and SCED problems. Currently, three strategies are employed to accelerate the computation. First is the method of identifying invalid constraints to reduce the scale of security constraints, as demonstrated in \cite{zhai2010fast, ma2020efficient}. Secondly, the security constraint can be treated as a type of lazy constraint \cite{lazy2021IBM}. This means that initially, these constraints are relaxed and not applied during the optimization process. They are only enforced when a potential solution is generated that violates the security constraint. At that point, the constraint/cutting plane is added to the problem, and the optimization process is repeated until a feasible solution that satisfies all constraints is found\cite{fu2007fast}. Nowadays, the method of machine learning is also used to accelerate calculations, as illustrated in reference \cite{xavier2021learning}.

The security constraint of DC power flow was originally expressed in the form of B-$\theta$. In order to simplify the modeling difficulty and reduce the number of variables and constraints, more literature is inclined to use Power Transfer Distribution Factors (PTDF) to represent the security constraint \cite{sahraei2016computationally, li2019enhanced}. PTDF is a linear displacement factor that assesses variations in the active power flow on a line as a result of variations in nodal power injection. On the MISO electricity market, the security constraints of RT-SCED are represented based on PTDF \cite{ma2009midwest}. Computational comparisons of the two formulations have been conducted in some published papers, such as \cite{sahraei2016computationally, zhang2021computational, hinojosa2017computational}. It is believed that with the number of decision variables and constraints decreasing in the PTDF-based formulation, the optimization speed will be much faster.  

Since the 1984 publication of Karmarkar's projective algorithm [11], interior point methods (IPMs) have been intensively investigated. Primal-dual IPMs are considered to be one of the most successful practical techniques \cite{lustig1994interior}. In general, simplex methods are typically more efficient for small and medium-sized problems, but interior point methods are typically more efficient for large-scale situations \cite{nocedal1999numerical}. The interior point method's overall performance \cite{kardovs2020structure} is mostly dependent on the sparse linear algebra solver. PTDF-based formulation generates a set of dense constraints and dealing with these denser constraints in the IPM typically necessitates a substantial amount of additional effort, making it difficult to solve the security-constrained optimization problem.

To the best of the authors' knowledge, there is currently no literature to exploit the sparsity structure of PTDF-based formulations. In PTDF-based formulations, only active powers are integrated into security constraints, while pure voltage phase angle formulations incorporate only voltage phase angles into security constraints. The regular B-$\theta$ formulation, referred to as the mixed formulation in this paper, includes both active powers and voltage phase angles in security constraints simultaneously.

In particular, the main contributions of this paper are as follows.
\begin{itemize}
  \item [1)] 
  As far as the authors' knowledge, the sparsity of the network security constraints in the primal-dual interior point technique is exploited for the first time in this paper. In IPM iterations, we further employ a matrix transformation technique to restore the sparsity of the linear system. 
  \item [2)]
  It is proved that the transformation approach is equivalent to solving the original optimization problem stated in terms of pure voltage angle, which maintains the sparsity structure but adds additional variables and constraints. Using such a voltage angle formulation, we can concentrate on the formulation rather than delving into the detailed realization of the solver. 
  \item [3)]
  Exhaustive numerical cases have been performed. results demonstrate that the solution speed of the mixed formulation is similar to that of the pure voltage angle formulation, but the numerical stability is much better. Additionally, the mixed formulation is faster than the PTDF formulation for optimal power flow (OPF) and security-constrained economic dispatch (SCED) for both the original problems and most of their outer approximations.
\end{itemize}

\section{Foundation}

\subsection{PTDF-based security constraints}
The PTDF-based security constraints for DC networks typically look like
\begin{equation}
  \label{eqn:ptdf_cons_origin}
  - \underline {{f_k}}  \le \sum\limits_{g \in \cal G} {{G_{kg}}{p_g}}  - \sum\limits_{d \in \cal D} {{D_{kd}}{p_d}}  \le \overline {{f_k}} ,\forall k \in \cal L
\end{equation}
where the forward and reverse maximum transmission power limits for the $k$th branch are denoted by $\overline{f_k}$ and $\underline{f_k}$, respectively. The injection shift factors for the $g$th generator and the $d$th load are represented by $G_{kg}$ and $D_{kd}$. The sets of branches, generators, and loads are denoted by $\mathcal{L}$, $\mathcal{G}$, and $\mathcal{D}$, respectively.

Equation (\ref{eqn:ptdf_cons_origin}) can be expressed in a more concise form as:
\begin{equation}
  \label{eqn:ptdf_cons_compact}
  - \underline {f}  \le PTDF \cdot p_{bus}  \le \overline {f} 
\end{equation}
where $\overline {f}$ and $\underline {f}$ are vectors that represent the transmission power limits for $n_l$ branches. $p_{bus}$ is the active power vector for $n_b$ buses, and the $PTDF$ matrix stores the injection shift factors for each bus to each branch, with a shape of $(n_l, n_b)$.

Formally, the $PTDF$ matrix can be formed with the incidence matrix $C$ and a diagonal branch reactance matrix $X=diag \{ x_k, \forall k \in \cal L \}$

\begin{equation}
  \label{eqn:ptdf_define}
  PTDF = B_f B_{bus}^{-1}  
\end{equation}
where $B_f := X^{-1} C$ is the coefficient from the voltage angle to branch power flow. $B_{bus} := C^T X^{-1} C$ is the nodal admittance matrix. To make a full-rank $B_{bus}$, the matrices $C$ and $X$ are built by eliminating a reference node. Note that both $B_f$ and $B_{bus}$ are sparse matrices. However, PTDF is highly dense because of the inverse of $B_{bus}$.

\subsection{Primal-dual interior point method for convex QPs}
For the following quadratic programming problem 
\begin{equation}
  \label{eqn:QP}
  \begin{array}{l}
    \min \frac{1}{2}{x^T}Gx + {c^T}x\\
    s.t.Ax \ge b
  \end{array}
\end{equation}

If the matrix $G$ is positive semi-definite, the optimization problem is referred to as a convex quadratic programming problem. By incorporating a slack variable $y$ and a dual variable $\lambda$, the Karush–Kuhn–Tucker (KKT) conditions can be derived as follows:
\begin{subequations}
  \label{eqn:KKT}
  \begin{align}
    Gx - {A^T}\lambda  + c = 0 \label{eqn:optimal_cond} \\ 
    Ax - y - b = 0\\
    {y_i}{\lambda _i} = 0, \forall i \label{eqn:element_wise}\\
    \lambda  \ge 0,y \ge 0 \label{eqn:inequality}
  \end{align}
\end{subequations}
where (\ref{eqn:element_wise}) is element-wise product.

The objective is to find solutions $(x^*, y^*, \lambda^*)$ such that KKT conditions (\ref{eqn:KKT}) hold. To achieve this, Primal-dual IPM applies a variant of Newton's method to (\ref{eqn:optimal_cond}) - (\ref{eqn:element_wise}) and modifies the search direction and step size at each iteration to ensure that the inequality (\ref{eqn:inequality}) is strictly satisfied.

Given an iteration point $(x_0, y_0, \lambda_0)$, the next pure Newton direction $(\Delta x, \Delta  y, \Delta \lambda)$ by solving the linear system as follows:
\begin{equation}
  \label{eqn:linear_system}
  \left[ {\begin{array}{*{20}{c}}
    G&{ - {A^T}}&0\\
    A&0&{ - I}\\
    0&Y&\Lambda 
    \end{array}} \right]\left[ {\begin{array}{*{20}{c}}
    {\Delta x}\\
    {\Delta \lambda }\\
    {\Delta y}
    \end{array}} \right] = \left[ {\begin{array}{*{20}{c}}
    { - {r_d}}\\
    { - {r_p}}\\
    { - {r_{\lambda y}}}
    \end{array}} \right]
\end{equation}
where dual residual ${r_d} = Gx_0 - {A^T}\lambda_0  + c $, primal residual ${r_p} = Ax_0 - y_0 - b $,  ${r_{\lambda y}} = \Lambda Ye$, $\Lambda  = diag(\lambda_0 )$, $Y = diag(y_0)$, $e = {[1,...,1]^T}$.

After eliminating $\Delta y$, the augment system form is 
\begin{subequations}
  \label{eqn:augment_system}
  \begin{align}
  \left[ {\begin{array}{*{20}{c}}
    G&{{A^T}}\\
    A&{ - {\Lambda ^{ - 1}}Y}
    \end{array}} \right]\left[ {\begin{array}{*{20}{c}}
    {\Delta x}\\
    { - \Delta \lambda }
    \end{array}} \right] = \left[ {\begin{array}{*{20}{c}}
    { - {r_d}}\\
    { - {r_p} - {\Lambda ^{ - 1}}{r_{\lambda y}}}
    \end{array}} \right] \\
    \Delta y = A\Delta x + {r_p}
  \end{align}
\end{subequations}

Since the coefficient matrix of the linear system is both sparse and symmetric, it can be efficiently solved via sparse LU decomposition. Furthermore, if the matrix G is similarly sparse, we can employ Gaussian elimination to generate the subsequent system \ref{eqn:cholesky} and solve it even more efficiently.
\begin{subequations}
  \label{eqn:cholesky}
  \begin{align}
    (G + {A^T}{Y^{ - 1}}\Lambda A)\Delta x =  - {r_d} - {A^T}{Y^{ - 1}}\Lambda ({r_p} + \Lambda {r_{\lambda y}}) \\
    \Delta y = A\Delta x + {r_p} \\
    \Delta \lambda  =  - {Y^{ - 1}}({r_{\lambda y}} + \Lambda \Delta y)
  \end{align}
\end{subequations}

The coefficient $G + {A^T}{Y^{ - 1}}\Lambda A$ is both sparse and positive semi-definite. This linear system can be solved using a modified sparse Cholesky decomposition.

In certain practical algorithms, such as Mehrotra's predictor-corrector method \cite{mehrotra1992implementation}, the predictor, corrector, and centering contributions are considered to improve the convergence speed of the IPM. The modification is restricted to the right side of the linear system and has no impact on the factorization of the coefficient matrix.

\section{PTDF-based security-constrained problem in IPM}

\subsection{Apply IPM to PTDF-based security-constrained problem}

A single-period PTDF-based SCED and OPF problem can be written as a QP problem
\begin{equation}
  \label{eqn:SCED}
  \begin{array}{l}
    \min \quad \frac{1}{2}\left[ {\begin{array}{*{20}{c}}
    {{p_{bus}^T}}&{{x_{aux}^T}}
    \end{array}} \right]\left[ {\begin{array}{*{20}{c}}
    {{G_{pp}}}&{G_{xp}^T}\\
    {{G_{xp}}}&{{G_{xx}}}
    \end{array}} \right]\left[ {\begin{array}{*{20}{c}}
    {{p_{bus}}}\\
    {{x_{aux}}}
    \end{array}} \right] \\ 
    + c_1^T{p_{bus}} + c_2^T{x_{aux}}\\
    s.t. \quad \left[ {\begin{array}{*{20}{c}}
    {{B_f}B_{bus}^{ - 1}}&J\\
    { - {B_f}B_{bus}^{ - 1}}&K\\
    M&L
    \end{array}} \right]\left[ {\begin{array}{*{20}{c}}
    {{p_{bus}}}\\
    {{x_{aux}}}
    \end{array}} \right] \ge \left[ {\begin{array}{*{20}{c}}
    {{b_1}}\\
    {{b_2}}\\
    {{b_3}}
    \end{array}} \right]
    \end{array}
\end{equation}
where $x_{aux}$ includes all decision variables except for bus active power $p_{bus}$. $b_1, b_2, b_3, c_1, c_2$ are constants vector and $G_{pp}, G_{xx}, G_{xp}, J, K, L, M$ are sparse constant matrices.

The KKT conditions of (\ref{eqn:SCED}) can be written as (\ref{eqn:SCED_KKT})
\begin{equation}
  \label{eqn:SCED_KKT}
  \begin{array}{l}
    \left[ {\begin{array}{*{20}{c}}
    {{G_{pp}}}&{G_{xp}^T}\\
    {{G_{xp}}}&{{G_{xx}}}
    \end{array}} \right]\left[ {\begin{array}{*{20}{c}}
    {{p_{bus}}}\\
    {{x_{aux}}}
    \end{array}} \right] - \\ \left[ {\begin{array}{*{20}{c}}
      {{B_f}B_{bus}^{ - 1}}&J\\
      { - {B_f}B_{bus}^{ - 1}}&K\\
      M&L
      \end{array}} \right]^T \left[ {\begin{array}{*{20}{c}}
    {{\lambda _1}}\\
    {{\lambda _2}}\\
    {{\lambda _3}}
    \end{array}} \right] 
    + \left[ {\begin{array}{*{20}{c}}
    {{c_1}}\\
    {{c_2}}
    \end{array}} \right] = 0\\
    \left[ {\begin{array}{*{20}{c}}
    {{B_f}B_{bus}^{ - 1}}&J\\
    { - {B_f}B_{bus}^{ - 1}}&K\\
    M&L
    \end{array}} \right]\left[ {\begin{array}{*{20}{c}}
    {{p_{bus}}}\\
    {{x_{aux}}}
    \end{array}} \right] - \left[ {\begin{array}{*{20}{c}}
    {{y_1}}\\
    {{y_2}}\\
    {{y_3}}
    \end{array}} \right] - \left[ {\begin{array}{*{20}{c}}
    {{b_1}}\\
    {{b_2}}\\
    {{b_3}}
    \end{array}} \right] = 0\\
    Y\Lambda  = 0\\
    \lambda  \ge 0,y \ge 0
    \end{array}
\end{equation}

The pure Newton direction can be derived by solving the linear equations following Gaussian elimination (\ref{eqn:SCED_Newton}).  
\begin{figure*}[hbtp] 
  \centering
  \begin{subequations}
    \label{eqn:SCED_Newton}
    \begin{align}
      \left[ {\begin{array}{*{20}{c}}
        {\left[ {\begin{array}{*{20}{c}}
        {{G_{pp}}}&{G_{xp}^T}\\
        {{G_{xp}}}&{{G_{xx}}}
        \end{array}} \right]}&{ - \left[ {\begin{array}{*{20}{c}}
        {B_{bus}^{ - T}B_f^T}&{ - B_{bus}^{ - T}B_f^T}&{{M^T}}\\
        {{J^T}}&{{K^T}}&{{L^T}}
        \end{array}} \right]}\\
        {\left[ {\begin{array}{*{20}{c}}
        {{B_f}B_{bus}^{ - 1}}&J\\
        { - {B_f}B_{bus}^{ - 1}}&K\\
        M&L
        \end{array}} \right]}&{{\Lambda ^{ - 1}}Y}
        \end{array}} \right]\left[ {\begin{array}{*{20}{c}}
        {\left[ {\begin{array}{*{20}{c}}
        {\Delta {p_{bus}}}\\
        {\Delta {x_{aux}}}
        \end{array}} \right]}\\
        {\left[ {\begin{array}{*{20}{c}}
        {\Delta {\lambda _1}}\\
        {\Delta {\lambda _2}}\\
        {\Delta {\lambda _3}}
        \end{array}} \right]}
        \end{array}} \right] = \left[ {\begin{array}{*{20}{c}}
        {\left[ {\begin{array}{*{20}{c}}
        { - {r_{d1}}}\\
        { - {r_{d2}}}
        \end{array}} \right]}\\
        {\left[ {\begin{array}{*{20}{c}}
        { - {r_{p1}}}\\
        { - {r_{p2}}}\\
        { - {r_{p3}}}
        \end{array}} \right] + {\Lambda ^{ - 1}}\left[ {\begin{array}{*{20}{c}}
        { - {r_{\lambda y1}}}\\
        { - {r_{\lambda y2}}}\\
        { - {r_{\lambda y3}}}
        \end{array}} \right]}
        \end{array}} \right] \label{eqn:dense_matrix} \\ 
        \left[ {\begin{array}{*{20}{c}}
          {\Delta {y_1}}\\
          {\Delta {y_2}}\\
          {\Delta {y_3}}
          \end{array}} \right] = \left[ {\begin{array}{*{20}{c}}
          {{B_f}B_{bus}^{ - 1}}&J\\
          { - {B_f}B_{bus}^{ - 1}}&K\\
          M&L
          \end{array}} \right]\left[ {\begin{array}{*{20}{c}}
          {\Delta {p_{bus}}}\\
          {\Delta {x_{aux}}}
          \end{array}} \right] + \left[ {\begin{array}{*{20}{c}}
          {{r_{p1}}}\\
          {{r_{p2}}}\\
          {{r_{p3}}}
          \end{array}} \right] \label{eqn:other_variable}
    \end{align}
  \end{subequations}
\end{figure*}

The constant parameters of (\ref{eqn:SCED_Newton}) can be calculated as (\ref{eqn:SCED_const}):
\begin{equation}
  \label{eqn:SCED_const}
  \begin{aligned}
    \left[ {\begin{array}{*{20}{c}}
    {{r_{d1}}}\\
    {{r_{d2}}}
    \end{array}} \right] =& \left[ {\begin{array}{*{20}{c}}
    {{G_{pp}}}&{G_{xp}^T}\\
    {{G_{xp}}}&{{G_{xx}}}
    \end{array}} \right]\left[ {\begin{array}{*{20}{c}}
    {{p_{bus}}}\\
    {{x_{aux}}}
    \end{array}} \right] - \\ 
    & \left[ {\begin{array}{*{20}{c}}
      {{B_f}B_{bus}^{ - 1}}&J\\
      { - {B_f}B_{bus}^{ - 1}}&K\\
      M&L
      \end{array}} \right]^T \left[ {\begin{array}{*{20}{c}}
    {{\lambda _1}}\\
    {{\lambda _2}}\\
    {{\lambda _3}}
    \end{array}} \right] + \left[ {\begin{array}{*{20}{c}}
    {{c_1}}\\
    {{c_2}}
    \end{array}} \right]\\
    \left[ {\begin{array}{*{20}{c}}
    {{r_{p1}}}\\
    {{r_{p2}}}\\
    {{r_{p3}}}
    \end{array}} \right] =& \left[ {\begin{array}{*{20}{c}}
    {{B_f}B_{bus}^{ - 1}}&J\\
    { - {B_f}B_{bus}^{ - 1}}&K\\
    M&L
    \end{array}} \right]\left[ {\begin{array}{*{20}{c}}
    {{p_{bus}}}\\
    {{x_{aux}}}
    \end{array}} \right] - \left[ {\begin{array}{*{20}{c}}
    {{y_1}}\\
    {{y_2}}\\
    {{y_3}}
    \end{array}} \right] - \left[ {\begin{array}{*{20}{c}}
    {{b_1}}\\
    {{b_2}}\\
    {{b_3}}
    \end{array}} \right]\\
    \left[ {\begin{array}{*{20}{c}}
    {{r_{\lambda y1}}}\\
    {{r_{\lambda y2}}}\\
    {{r_{\lambda y3}}}
    \end{array}} \right] =& \Lambda Ye\\
    \Lambda  = diag(\lambda )&,Y = diag(y),e = {[1,...,1]^T}
    \end{aligned}
\end{equation}

The calculation of ${B_f} B_{bus}^{-1}$ results in a highly dense matrix, which poses challenges when using traditional sparse decomposition methods to solve the resulting linear system. As factorization is the primary task in each iteration of the IPM, this can make the algorithm inefficient.

\subsection{Matrix transformation to restore sparsity}

In order to restore the sparsity of (\ref{eqn:dense_matrix}), we find a matrix transformation method to achieve this purpose. Here is an example to provide a hint: Consider the following linear equation (\ref{eqn:example}).
\begin{equation}
  \label{eqn:example}
  \left[ {\begin{array}{*{20}{c}}
    G&{{B_{bus}^{ - T}}{B_f^T}}\\
    {B_f{B_{bus}^{ - 1}}}&0
    \end{array}} \right]\left[ {\begin{array}{*{20}{c}}
    {\Delta x}\\
    {\Delta \lambda }
    \end{array}} \right] = \left[ {\begin{array}{*{20}{c}}
    a\\
    b
    \end{array}} \right]
\end{equation}

It is equivalent to solving (\ref{eqn:example_equ}).
\begin{equation}
  \label{eqn:example_equ}
  \left[ {\begin{array}{*{20}{c}}
    {{B_{bus}^T}GB_{bus}}&{{B_f^T}}\\
    B_f&0
    \end{array}} \right]\left[ {\begin{array}{*{20}{c}}
    {{B_{bus}^{ - 1}}\Delta x}\\
    {\Delta \lambda }
    \end{array}} \right] = \left[ {\begin{array}{*{20}{c}}
    {{B_{bus}^T}a}\\
    b
    \end{array}} \right]
\end{equation}

The coefficient matrix of (\ref{eqn:example_equ}) is sparse because of the sparsity of ${B_{bus}^T}GB_{bus}$. The sparsity is recovered through this matrix transformation.

By applying the matrix transformation method to (\ref{eqn:dense_matrix}), we can obtain an equivalent expression of (\ref{eqn:SCED_Newton}) as shown in (\ref{eqn:transformation}), where the sparsity of the coefficient matrix is recovered due to the sparsity of ${B_{bus}^T}GB_{bus}$. As ${B_{bus}^T}GB_{bus}$ is both sparse and positive semi-definite, we can efficiently solve the resulting linear system using modified sparse Cholesky decomposition.

\begin{figure*}[hbtp]
  \begin{subequations}
    \label{eqn:transformation}
    \begin{align}
      \left[ {\begin{array}{*{20}{c}}
        {B_{bus}^T{G_{pp}}{B_{bus}}}&{B_{bus}^TG_{xp}^T}&{B_f^T}&{ - B_f^T}&{B_{bus}^T{M^T}}\\
        {{G_{xp}}{B_{bus}}}&{{G_{xx}}}&{{J^T}}&{{K^T}}&{{L^T}}\\
        {{B_f}}&J&{ - \Lambda _1^{ - 1}{Y_1}}&0&0\\
        { - {B_f}}&K&0&{ - \Lambda _2^{ - 1}{Y_2}}&0\\
        {M{B_{bus}}}&L&0&0&{ - \Lambda _3^{ - 1}{Y_3}}
        \end{array}} \right]\left[ {\begin{array}{*{20}{c}}
        {\Delta {\theta}}\\
        {\Delta {x_{aux}}}\\
        { - \Delta {\lambda _1}}\\
        { - \Delta {\lambda _2}}\\
        { - \Delta {\lambda _3}}
        \end{array}} \right] = \left[ {\begin{array}{*{20}{c}}
        { - B_{bus}^T{r_{d1}}}\\
        { - {r_{d2}}}\\
        { - {r_{p1}} - \Lambda _1^{ - 1}{r_{\lambda y1}}}\\
        { - {r_{p2}} - \Lambda _2^{ - 1}{r_{\lambda y2}}}\\
        { - {r_{p3}} - \Lambda _3^{ - 1}{r_{\lambda y3}}}
        \end{array}} \right] \\
        \Delta {p_{bus}} = {B_{bus}}\Delta \theta  \label{eqn:mapping} \\
        (\ref{eqn:other_variable})
      \end{align}
    \end{subequations}

\end{figure*}

\subsection{Security constraints formulations}

It is a challenging and impractical task to modify the IPM iterations in a well-established solver. Moreover, SCED is a multi-stage problem, which further complicates the formulation of the transformed equation (\ref{eqn:transformation}) due to its more complex structure.

Note that equation (\ref{eqn:mapping}) represents a linear mapping of the gradient of $p_{bus}$. By performing the same linear mapping of $p_{bus}$ to the original problem (\ref{eqn:SCED}), we can obtain a new formulation (\ref{eqn:pure_theta}) that is defined in terms of pure voltage angle $\theta$ instead of nodal injection $p_{bus}$. It is straightforward to verify that the Newton step of (\ref{eqn:pure_theta}) is also given by equation (\ref{eqn:transformation}).

\begin{equation}
  \label{eqn:pure_theta}
  \begin{array}{l}
    \min \frac{1}{2}\left[ {\begin{array}{*{20}{c}}
    \theta^T &{{x_{aux}^T}}
    \end{array}} \right]\left[ {\begin{array}{*{20}{c}}
    {B_{bus}^T{G_{pp}}{B_{bus}}}&{B_{bus}^TG_{xp}^T}\\
    {{G_{xp}}{B_{bus}}}&{{G_{xx}}}
    \end{array}} \right] \left[ {\begin{array}{*{20}{c}}
    \theta \\
    {{x_{aux}}}
    \end{array}} \right] \\
    + c_1^T{B_{bus}}\theta  + c_2^T{x_{aux}}\\
    s.t. \left[ {\begin{array}{*{20}{c}}
    {{B_f}}&J\\
    { - {B_f}}&K\\
    {M{B_{bus}}}&L
    \end{array}} \right]\left[ {\begin{array}{*{20}{c}}
    \theta \\
    {{x_{aux}}}
    \end{array}} \right] \ge \left[ {\begin{array}{*{20}{c}}
    {{b_1}}\\
    {{b_2}}\\
    {{b_3}}
    \end{array}} \right]
    \end{array}
\end{equation}

By transforming all nodal injection variables to voltage angles, the PTDF-based SCED and OPF problems can be solved. However, in most SCED and OPF problems, the constraints (such as upper and lower power constraints of generators) and objective functions are often about the generator power $p_g$ rather than the nodal net injection $p_{bus}$, as shown in (\ref{eqn:PTDF_formulation}). 

\textbf{PTDF Formulation:}
\begin{equation}
  \label{eqn:PTDF_formulation}
  \begin{aligned}
    \text{min} \quad & \sum_{g \in \cal G} \left( a_g {p_g}^2 + b_g p_g + c_g \right) \\
    \text{s.t} \quad & \underline{p_g} \leq p_g \leq \overline{p_g}, \quad \forall g \in \cal G \\
    & \text{Security Constraints (\ref{eqn:ptdf_cons_origin}) }
  \end{aligned}
  \end{equation}
where $a_g$, $b_g$, and $c_g$ are the coefficients of generator cost function. The $\overline {p_g}$ and $\underline {p_g}$ are the generator power limits for the $g$th generator.

Two techniques are proposed for converting the standard PTDF formulation (\ref{eqn:PTDF_formulation}) into the form of (\ref{eqn:pure_theta}). In this paper, they are referred to as the pure voltage angle formulation and the mixed formulation, respectively.

The pure voltage angle formulation maintains the structure of (\ref{eqn:pure_theta}) exactly. In this formulation, the nodal injection power is the aggregate of the generator power and local load. The generator power constraint can be conveniently transformed into an equivalent nodal injection power constraint. Moreover, the convex cost functions of several generators can be consolidated into an equivalent convex nodal cost function, which can be linearized piecewise. Following the optimization process, the nodal net injection power can be decomposed into individual generators based on their marginal cost and power constraints.

When only one generator is connected to a bus, the nodal cost function can be expressed as a quadratic function of the net nodal injection $p_b$. The following formulation assumes at most one generator per bus to enable a fair comparison of the efficiency of quadratic programming. 

\textbf{Pure Voltage Angle Formulation:}
\begin{equation}
  \label{eqn:theta_formulation}
  \begin{aligned}
    \text{min} \quad & \sum_{b \in \cal B} \left[ a_b {(B_{bus, b} \theta)}^2 + b_b B_{bus, b} \theta + c_b \right] \\
    \text{s.t} \quad & \underline{p_b} \leq B_{bus, b} \theta \leq \overline{p_b}, \quad \forall b \in \cal B \\
    &  -\underline {f} \leq B_f \theta \leq  \overline{f}
  \end{aligned}
\end{equation}
where $B_{bus, b}$ represents the $b$th row of $B_{bus}$ matrix.If the generator power and load at the $b$th bus are $p_g$ and $p_d$ respectively, the cost coefficients will be $a_b = a_g$, $b_b = b_g + 2a_g p_d$, and $c_b = c_g + b_g p_d + a_g p_d^2$. The nodal power limits will be $\underline{p_b} = \underline{p_g} + p_d$ and $\overline{p_b} = \overline{p_g} + p_d$.

The mixed formulation preserves both voltage angle and generator output in the optimization problem. The voltage angle represents security constraints, while the generator powers indicate the remaining constraints and objectives. However, since the problem is not fully transformed into a voltage angle problem, the relationship between voltage angle and nodal net injection, $p_{bus} = B_{bus} \theta$, must be included as a constraint. This formulation is also referred to as the B-$\theta$ formulation.

\textbf{Mixed Formulation:}
\begin{equation}
  \label{eqn:mixed_formulation}
  \begin{aligned}
    \text{min} \quad & \sum_{g \in \cal G} \left( a_g {p_g}^2 + b_g p_g + c_g \right) \\
    \text{s.t} \quad & \underline{p_g} \leq p_g \leq \overline{p_g}, \quad \forall g \in \cal G \\
    &  -\underline {f} \leq B_f \theta \leq  \overline{f} \\
    & \sum_{g \in {\cal G}_b} p_g - \sum_{d \in {\cal D}_b} p_d = B_{bus,b} \theta, \quad \forall b \in \cal B 
  \end{aligned}
\end{equation}
where ${\cal G}_b$ and ${\cal D}_b$ are the set of generators and loads at bus $b$.

The comparison of a single-period OPF problem in the aforementioned formulations is presented in Table \ref{tab:compare}.

\begin{table}[htb]
  \renewcommand\arraystretch{1.25}
  \caption{Comparison of three formulations for single-period OPF problem}
  \begin{threeparttable}
    \begin{tabular}{lccc}
    \hline
    \hline
    &
      \textbf{\begin{tabular}[c]{@{}c@{}}PTDF\\ formulation\end{tabular}} &
      \textbf{\begin{tabular}[c]{@{}c@{}}Pure voltage\\ angle formulation\end{tabular}} &
      \textbf{\begin{tabular}[c]{@{}c@{}}Mixed \\ formulation\end{tabular}} \\ \hline
    \textbf{Variables}   & $\left|\cal{G}\right|$ & $\left|\cal{B}\right|$   & $\left|\cal{G}\right|+\left|\cal{B}\right|$ \\
    \textbf{Constraints} & $2\left|\cal{L}\right| + 2\left|\cal{G}\right|$ & $2\left|\cal{L}\right| + 2\left|\cal{B}\right|$   & $2\left|\cal{L}\right|+ 2 \left|\cal{G}\right| + \left|\cal{B}\right|$ \\
    \textbf{Sparsity}    & $\times$   & $\surd$   & $\surd$        \\ 
    \hline
    \hline
    \end{tabular}
    \begin{tablenotes}   
      \footnotesize              
      \item[*]  Only count the number of variables and constraints associated with the security constraint and generation power limits constraints.
      \item[**]  A single equality constraint is counted as one constraint.
    \end{tablenotes} 
  \end{threeparttable} 
  \label{tab:compare}
\end{table}

\section{Numerical results}

This section presents the numerical results of the three formulations. First, we provide an overview of the coefficient matrices' structure for each of the three formulations. Next, we present the case studies for OPF and SCED problems. We also assess the formulations' ability to handle security constraints using outer approximations. Finally, we discuss the limitations of the proposed formulations in mixed-integer problems, such as the SCUC problem. All the source code used in this study is available on GitHub \cite{zhiyuan2023github}.

\subsection{Implementation details}

We model the cost of each generator as a convex quadratic function. The test cases come with built-in data, including the active load on each bus, the status of generators, the minimum and maximum generator output limits, and the power flow limits of each branch. For branches without power flow limit data, security constraints do not consider them. The detailed parameters for each case can be found at \cite{zhiyuan2023github}.

In the single-period OPF problem, the objective function aims to minimize the cost of the generator currently operating within a single period. Constraints include the limits on generator output power and branch power flow. In order to make the three formulations comparable, generators, which are located in the same bus, will be merged into one generator by setting the same type of parameters. The cost function w.r.t. nodal net injection can be written as a convex quadratic function.

In the SCED problem, we divide a day into 24 periods, each lasting one hour. The objective function is to minimize the total cost of generators that will be operated in the next 24 periods. In addition to the constraints in the OPF problem, the ramping constraints of each generator should also be met. The bus loads in each period are randomly set between 0.95 and 1.05 times the original loads. 

In the SCUC problem, we include the parameters used in SCED, as well as fixed startup and shutdown costs, a reserve ratio of 0.1 of total loads, and a minimum startup and shutdown duration of 4 hours. The SCUC formulation also adds additional constraints such as the reserve constraint and the minimum startup/shutdown duration constraint to the SCED problem.

The cases used in this study are all built-in test cases of MATPOWER \cite{zimmerman2010matpower} with more than 1000 buses, which include the ACTIV Synthetic Grid Test Cases \cite{birchfield2016grid}, the Polish System Test Cases, the PEGASE European System Test Cases \cite{josz2016ac, fliscounakis2013contingency}, and the RTE French System Test Cases\cite{josz2016ac}.

The MATPOWER functions \textit{makePTDF} and \textit{makeBdc} are used to calculate $PTDF$, $B_{bus}$ and $B_{f}$. SCUC is a MIQP problem while OPF and SCED are convex QP problems. These problems are coded with YALMIP (version: R20200930) in MATLAB R2019a, solved by Gurobi 9.0.3, and executed on an Intel i7-9700 processor with 8 threads and 32GB of memory.

\subsection{Sparsity structure of coefficient matrices}

This section discusses how the PTDF-based formulation impacts the sparsity of the coefficient matrix and the efficiency of the interior point method.

\begin{figure}[htbp!]
  \centering
  \includegraphics[width=0.45\textwidth]{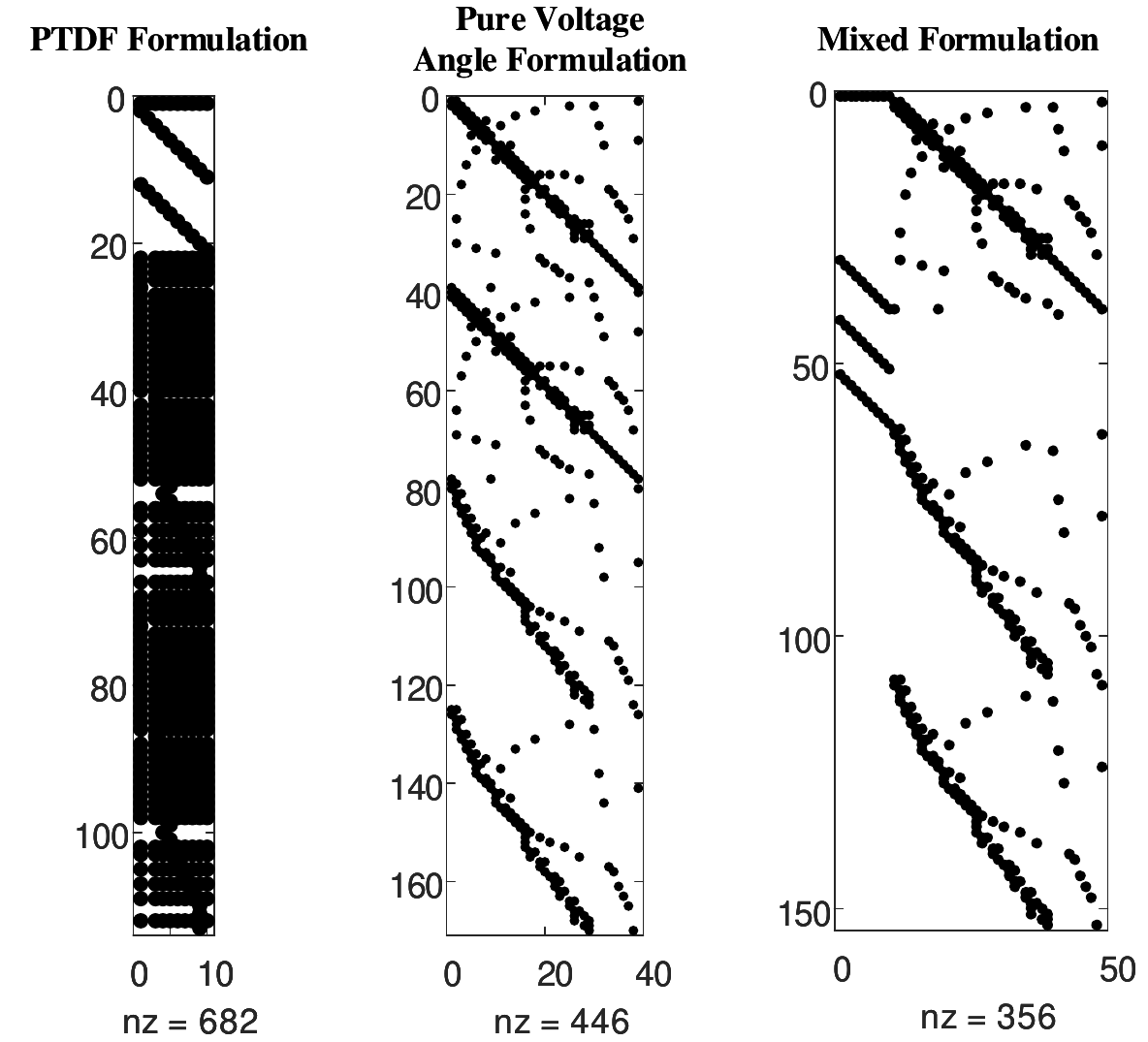}
  \caption{Illustration of the non-zero elements distribution in the coefficient matrices of three formulations used for the OPF problem in the IEEE 39-bus power system. These matrices were obtained prior to the presolve process.}
  \label{fig:structure}
\end{figure}

To illustrate this, we analyze a simple IEEE 39-bus power system with 39 buses, 10 generators, and 46 transmission lines. Figure \ref{fig:structure} displays the coefficient matrix structure of three different formulations for the single-period OPF problem.

In the PTDF formulation, the number of variables is $\left|\cal{G}\right|= 10$ , and the number of constraints is $2\left|\cal{L}\right|+2\left|\cal{G}\right|=112$. Due to the high density of the PTDF matrix, the number of non-zero elements is 682. In pure voltage angle formulation, the number of variables is $\left|\cal{B}\right|= 39$ , the number of constraints is $2\left|\cal{L}\right|+2\left|\cal{B}\right|=170$, and the number of non-zero elements is 446. In mixed formulation, the number of variables is $\left|\cal{G}\right|+\left|\cal{B}\right|= 49$ , the number of constraints is $2\left|\cal{L}\right|+ 2 \left|\cal{G}\right| + \left|\cal{B}\right|=151$, and the number of non-zero elements is 356. It is evident that $PTDF < Pure \approx Mixed$ in terms of coefficient matrix size, whereas $PTDF > Pure > Mixed$ holds true in terms of coefficient matrix density.

In order to confirm that the above conclusion holds true in the larger test cases and that the dense coefficient matrices would slow down IPM, we evaluate the performance of solving a single-period OPF in Table \ref{tab:OPF_result}'s five large test cases.

\begin{table*}[htbp]
  \renewcommand\arraystretch{1.1}
  \caption{Comparison of size, density and solving time of three formulations for solving the single-period OPF problem}
  \begin{threeparttable}
    \begin{tabular}{cccccccccc}
      \hline\hline
       &  & \multicolumn{8}{c}{PTDF Based Formulation} \\ \cline{3-10} 
       &  & Constraints & Variables & A NZ & A NZ (\%) & AA' NZ & Factor Ops & Solver Time(s) & Barrier Time(s) \\ \hline
      \multirow{5}{*}{\begin{tabular}[c]{@{}c@{}}Without \\ Presolve\end{tabular}} & case1951rte & \textbf{4231} & \textbf{358} & 9.23E+05 & 60.92\% & 9.23E+05 & 3.63E+08 & 0.77 & 0.49 \\
       & case6470rte & \textbf{6989} & \textbf{634} & 3.07E+06 & 69.19\% & 3.07E+06 & 2.09E+09 & 2.35 & 1.56 \\
       & case6515rte & \textbf{6893} & \textbf{563} & 2.72E+06 & 70.11\% & 2.72E+06 & 1.64E+09 & 2.29 & 1.57 \\
       & case\_ACTIVSg10k & \textbf{22231} & \textbf{1455} & 2.46E+07 & 76.18\% & 2.46E+07 & 3.74E+10 & 26.98 & 20.60 \\
       & case\_ACTIVSg25k & \textbf{49521} & \textbf{2753} & 1.07E+08 & 78.64\% & 1.07E+08 & 3.05E+11 & 554.55 & 522.92 \\ \hline
      \multirow{5}{*}{\begin{tabular}[c]{@{}c@{}}With\\ Presolve\end{tabular}} & case1951rte & \textbf{192} & \textbf{386} & 6.78E+04 & 91.50\% & 1.83E+04 & 2.38E+06 & 0.99 & 0.24 \\
       & case6470rte & \textbf{426} & \textbf{728} & 2.68E+05 & 86.43\% & 9.05E+04 & 2.59E+07 & 1.54 & 0.41 \\
       & case6515rte & \textbf{470} & \textbf{663} & 2.62E+05 & 84.20\% & 1.10E+05 & 3.47E+07 & 1.91 & 0.56 \\
       & case\_ACTIVSg10k & \textbf{769} & \textbf{783} & 5.36E+05 & 89.03\% & 2.95E+05 & 1.52E+08 & 8.24 & 2.03 \\
       & case\_ACTIVSg25k & \textbf{3746} & \textbf{3872} & 1.03E+07 & 71.08\% & 7.01E+06 & 1.75E+10 & 64.81 & 35.12 \\ \hline\hline
       &  & \multicolumn{8}{c}{Pure Voltage Angle Formulation} \\ \cline{3-10} 
       &  & Constraints & Variables & A NZ & A NZ (\%) & AA' NZ & Factor Ops & Solver Time(s) & Barrier Time(s) \\ \hline
      \multirow{5}{*}{\begin{tabular}[c]{@{}c@{}}Without \\ Presolve\end{tabular}} & case1951rte & 8100 & 1951 & 2.18E+04 & 0.14\% & 1.49E+05 & 1.56E+07 & 0.31 & 0.22 \\
       & case6470rte & 19160 & 6470 & 5.76E+04 & 0.05\% & 3.05E+05 & 4.72E+07 & 0.44 & \textbf{0.25} \\
       & case6515rte & 19292 & 6515 & 5.80E+04 & 0.05\% & 3.06E+05 & 4.88E+07 & 0.48 & 0.27 \\
       & case\_ACTIVSg10k & 40488 & 10000 & 1.10E+05 & 0.03\% & 6.25E+05 & 1.48E+08 & 1.55 & 0.91 \\
       & case\_ACTIVSg25k & 96660 & 25000 & 2.64E+05 & 0.01\% & 1.38E+06 & 5.07E+08 & 6.32 & 2.81 \\ \hline
      \multirow{5}{*}{\begin{tabular}[c]{@{}c@{}}With\\ Presolve\end{tabular}} & case1951rte & 1812 & 2170 & 6.58E+03 & 0.17\% & 1.79E+04 & 1.18E+06 & 0.31 & \textbf{0.17} \\
       & case6470rte & 3680 & 4314 & 1.51E+04 & 0.10\% & 4.09E+04 & 5.47E+06 & 0.47 & 0.21 \\
       & case6515rte & 3512 & 4075 & 1.45E+04 & 0.10\% & 3.98E+04 & 5.08E+06 & 0.39 & \textbf{0.20} \\
       & case\_ACTIVSg10k & 9948 & 10648 & 3.61E+04 & 0.03\% & 1.03E+05 & 2.16E+07 & 1.04 & 0.35 \\
       & case\_ACTIVSg25k & 22190 & 24943 & 8.35E+04 & 0.02\% & 2.15E+05 & 9.43E+07 & 3.55 & 0.74 \\ \hline\hline
       &  & \multicolumn{8}{c}{Mixed Formulation} \\ \cline{3-10} 
       &  & Constraints & Variables & A NZ & A NZ (\%) & AA' NZ & Factor Ops & Solver Time(s) & Barrier Time(s) \\ \hline
      \multirow{5}{*}{\begin{tabular}[c]{@{}c@{}}Without \\ Presolve\end{tabular}} & case1951rte & 6867 & 2309 & \textbf{1.65E+04} & 0.10\% & 8.11E+04 & 6.36E+06 & 0.18 & \textbf{0.12} \\
       & case6470rte & 13960 & 7104 & \textbf{3.76E+04} & 0.04\% & 1.43E+05 & 1.62E+07 & 0.33 & \textbf{0.25} \\
       & case6515rte & 13905 & 7078 & \textbf{3.75E+04} & 0.04\% & 1.43E+05 & 1.66E+07 & 0.34 & \textbf{0.26} \\
       & case\_ACTIVSg10k & 33400 & 11455 & \textbf{8.12E+04} & 0.02\% & 3.14E+05 & 6.49E+07 & 1.08 & \textbf{0.45} \\
       & case\_ACTIVSg25k & 77168 & 27753 & \textbf{1.90E+05} & 0.01\% & 6.95E+05 & 2.51E+08 & 2.31 & \textbf{1.48} \\ \hline
      \multirow{5}{*}{\begin{tabular}[c]{@{}c@{}}With\\ Presolve\end{tabular}} & case1951rte & 1080 & 1434 & \textbf{4.67E+03} & 0.30\% & 1.02E+04 & 1.01E+06 & 0.30 & \textbf{0.17} \\
       & case6470rte & 2402 & 3033 & \textbf{1.18E+04} & 0.16\% & 2.86E+04 & 4.19E+06 & 0.32 & \textbf{0.18} \\
       & case6515rte & 2391 & 2951 & \textbf{1.16E+04} & 0.16\% & 2.84E+04 & 4.16E+06 & 0.35 & \textbf{0.20} \\
       & case\_ACTIVSg10k & 7727 & 8425 & \textbf{3.01E+04} & 0.05\% & 7.17E+04 & 2.24E+07 & 0.49 & \textbf{0.33} \\
       & case\_ACTIVSg25k & 18293 & 21044 & \textbf{7.46E+04} & 0.02\% & 1.71E+05 & 8.86E+07 & 0.96 & \textbf{0.70} \\ \hline\hline
      \end{tabular}
  \end{threeparttable}
  \label{tab:OPF_result}
  \end{table*}

Table \ref{tab:OPF_result} presents the results of three different formulations, both with and without the presolve process. It is worth noting that the optimal value is identical for all three formulations. To compare the formulation sizes, the first two columns show the number of constraints and variables in each formulation, respectively. The third and fourth columns represent the density of the formulation, in terms of the number and proportion of non-zero coefficient matrix $A$ elements, respectively. Column 5 displays the number of non-zero matrix members that require decomposition in each IPM iteration. Specifically, Gurobi calculates the number of non-zero elements in the lower triangular matrix of $AA^T$, and column 6 indicates the number of floating-point operations required to factor it. Generally, a single core of current processors can execute approximately $5\times10^9$ floating-point operations per second, and most problems converge in no more than 50 IPM iterations, so the solution time can be estimated. The final two columns show the solver time, which includes both the barrier time and the data conversion time between Gurobi and YALMIP.

Table \ref{tab:OPF_result} reveals that the PTDF formulation has a much smaller problem size, a higher density, and a slower solving speed compared to the other two formulations. For instance, using the \textit{case\_ACTIVSg25k} as an example, the size of the PTDF formulation is approximately $1/20$ of that of the pure angle and mixed formulation. However, the number of non-zero elements in matrix $A$ is over 100 times greater, and the proportion of non-zero elements in the PTDF formulation is over 70\%, which is significantly denser than the 0.02\% of the other formulations. As a result, the factorization process in each IPM iteration for the PTDF formulation requires more than 100 times the number of floating-point operations, and the barrier time is 50 times longer than that of the other formulations.

The performance difference between the pure voltage angle formulation and the mixed formulation is insignificant. However, after presolve, the size and number of non-zero elements in the mixed formulation are considerably reduced, resulting in a faster solving speed. Since the pure voltage angle method requires the solver to confirm that the quadratic term matrix in the objective function (\ref{eqn:pure_theta}) is a positive semi-definite matrix, it can be solved using the convex QP algorithm. However, in some cases, the solver is unable to determine a positive semi-definite matrix and hence erroneously applies the non-convex QP technique. Therefore, we strongly believe that the mixed formulation is superior to the pure voltage angle formulation. Consequently, we will focus exclusively on the PTDF and mixed formulations in the following section.

\subsection{Performance comparison}

\begin{figure}[htbp!]
  \centering
  \includegraphics[width=0.4\textwidth]{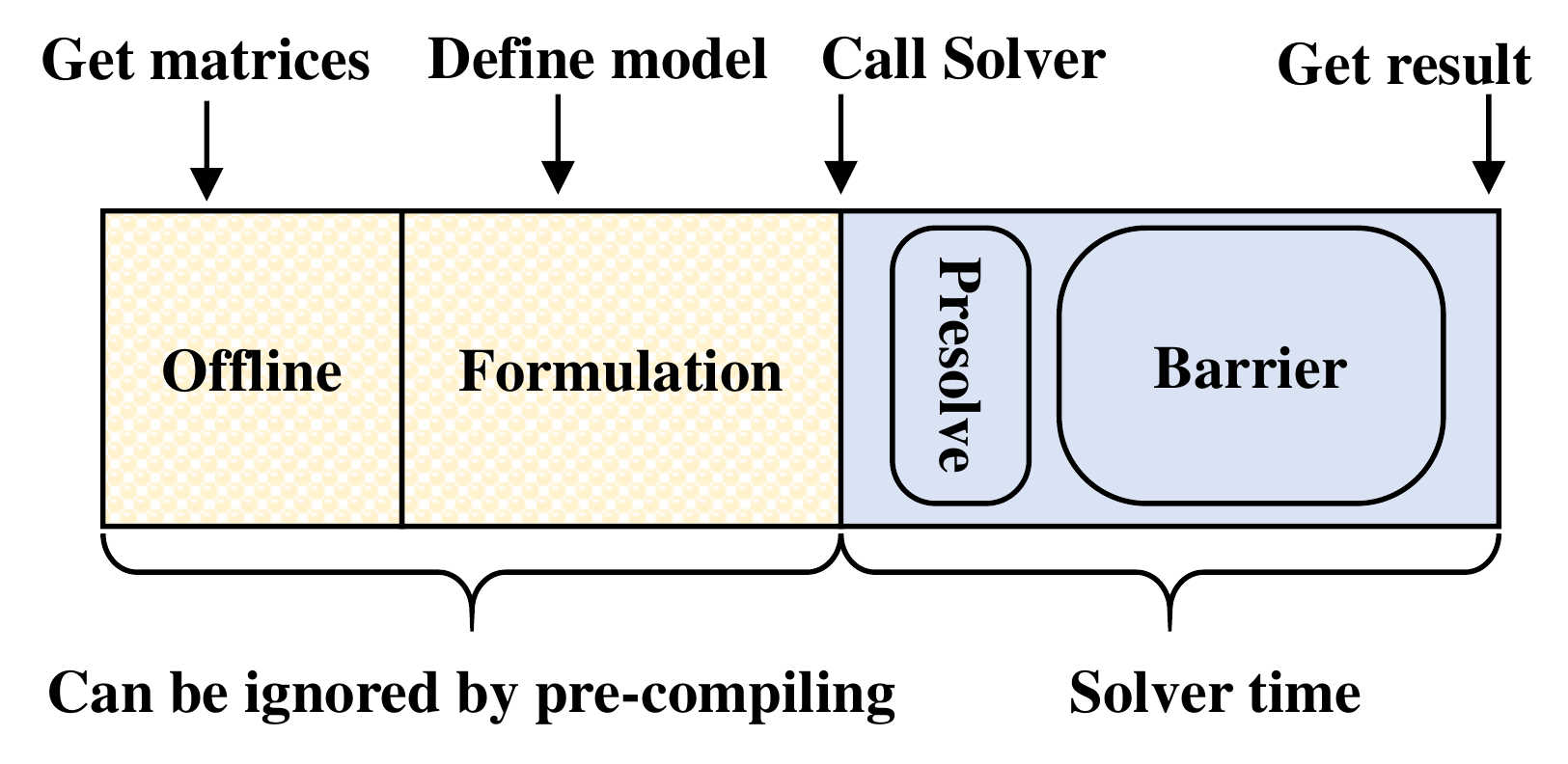}
  \caption{Time composition for solving optimization problems. The offline and formulation time can be ignored by pre-compiling. The main component of solver time is barrier time, which is crucial in measuring the efficacy of various formulations. }
  \label{fig:time}
\end{figure}

Table \ref{tab:SCED_result} displays the solution times for the PTDF and mixed formulations of the SCED problem. The detailed time composition for solving an optimization problem is illustrated in Fig. \ref{fig:time}. Offline time denotes the time needed to construct the $PTDF$, $B_{bus}$, and $B_f$ matrices. Formulation time is the time necessary to model an optimization problem using YALMIP. Solver time is the duration between calling the solver and receiving the final result. The barrier time represents the duration of performing IPM. Factor Ops is the number of floating-point operations required per iteration for matrix factorization in a single iteration. The model formulation time can be negligible by creating a pre-compiled model.

\begin{table*}[htbp]
  \renewcommand\arraystretch{1.1}
  \caption{Detailed time comparison of PTDF and mixed formulations for SCED problem}
  \begin{threeparttable}
    \begin{tabular}{cccccccccccc}
      \hline\hline
       & \multicolumn{5}{c}{PTDF Based Formulation} &  & \multicolumn{5}{c}{Mixed Formulation} \\ \cline{2-6} \cline{8-12} 
       & Offline(s) & Formulation(s) & Solver(s) & Barrier(s) & Factor Ops &  & Offline(s) & Formulation(s) & Solver(s) & Barrier(s) & Factor Ops \\ \hline
      case1354pegase & 0.30 & 7.09 & 8.44 & 5.31 & 4.26E+09 &  & 0.00 & 12.76 & 1.85 & \textbf{1.50} & 9.96E+08 \\
      case1951rte & 0.53 & 13.10 & 11.24 & 6.05 & 1.22E+10 &  & 0.00 & 25.50 & 2.07 & \textbf{1.73} & 2.32E+09 \\
      case3375wp & 1.82 & 28.71 & 28.47 & 15.14 & 1.90E+10 &  & 0.00 & 53.46 & 5.30 & \textbf{4.80} & 5.25E+09 \\
      case6470rte & 8.26 & 55.25 & 334.93 & 315.11 & 2.19E+12 &  & 0.00 & 190.12 & 7.36 & \textbf{6.68} & 1.59E+10 \\
      case6515rte & 8.35 & 47.57 & 331.84 & 314.74 & 2.42E+12 &  & 0.00 & 165.72 & 7.30 & \textbf{6.65} & 1.60E+10 \\ \hline\hline
      \end{tabular}
      \begin{tablenotes}   
        \footnotesize              
        \item[*]  All results are obtained with presolve process. 
      \end{tablenotes} 
  \end{threeparttable}
  \label{tab:SCED_result}
\end{table*}

\begin{figure}[htbp!]
  \centering
  \includegraphics[width=0.5\textwidth]{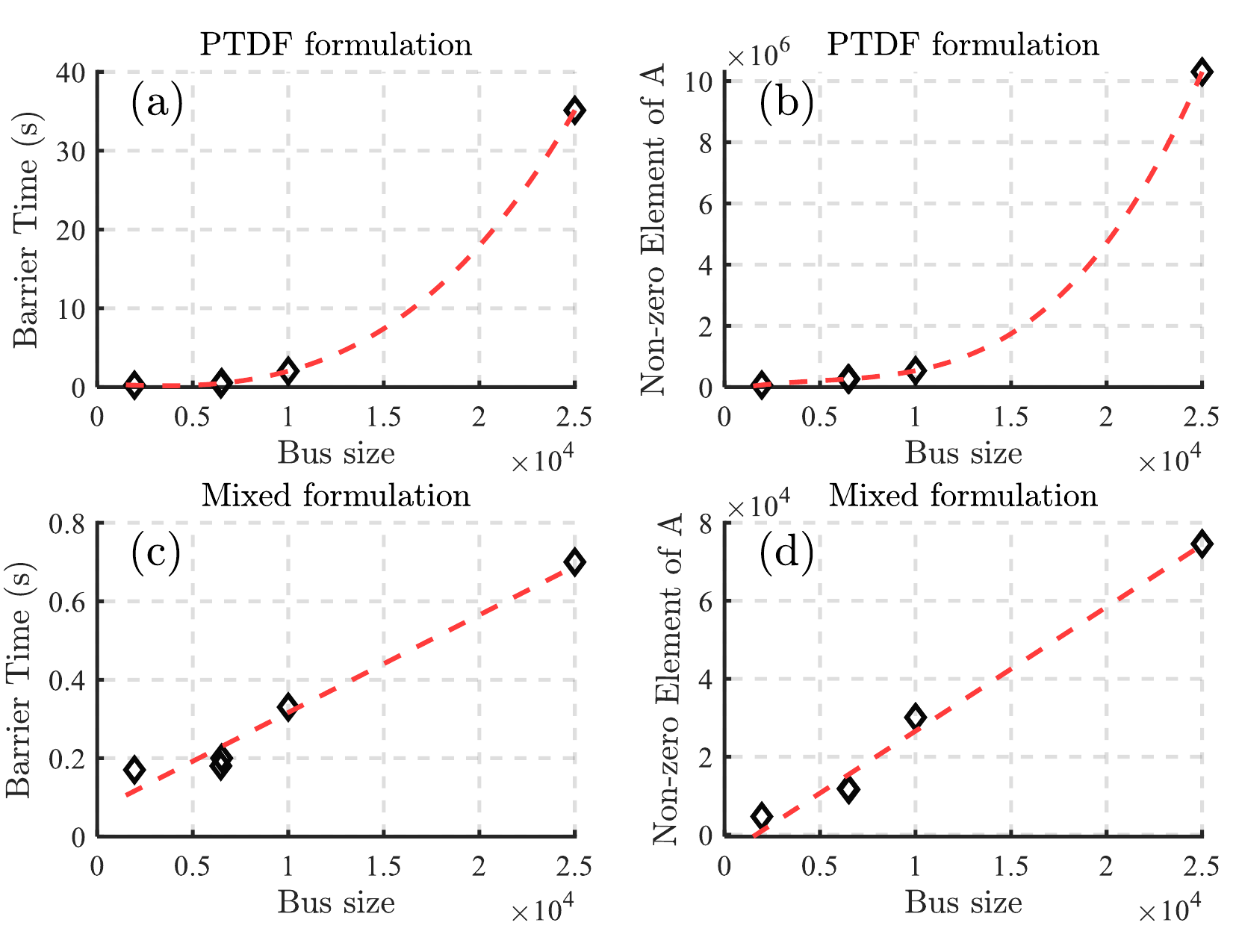}
  \caption{Barrier time and non-zero elements for varying the number of buses in PTDF and mixed formulation. The red dashed lines are the fitted curves.}
  \label{fig:fit}
\end{figure}

The results show that as the number of variables and constraints increases, the formulation time of the mixed formulation increases. However, this can be neglected by pre-compiling the model. On the other hand, the solution time is significantly reduced compared to the PTDF formulation. As depicted in Fig. \ref{fig:fit}, as the number of buses in the test cases for the OPF problem increases with presolve, both the barrier time and the number of non-zero elements in the coefficient matrix of the IPM algorithm increase. The barrier time and non-zero elements of PTDF increase superlinearly as the number of buses increases, but those of the mixed formulation increase approximately linearly. Consequently, the mixed formulation provides greater comparative benefits when dealing with problems on a larger scale.

\subsection{Performance on outer approximations of security constraints}

Many branch constraints in security constraints are redundant. Extensive research has been conducted on the technique for excluding redundant branch constraints. After removing redundant constraints, the resulting optimization problem becomes an outer approximation of the original problem. In this work, we adopt two strategies to obtain outer approximations of the SCED problem for both the PTDF and mixed formulations: 1) randomly removing uncongested branch constraints, and 2) sequentially removing the most uncongested branch constraints. The results are presented in Table \ref{tab:SCED_del_result}.

\begin{table*}[htbp]
  \renewcommand\arraystretch{1.1}
      \caption{Performances of PTDF and mixed formulation in the outer approximations of the SCED problem}
      \begin{threeparttable}
          \begin{tabular}{ccccccccccccccccc}
            \hline\hline
            \multicolumn{2}{c}{\multirow{3}{*}{\diagbox{Case}{Time(s)}{Del. Ratio}}} & \multicolumn{7}{c}{\multirow{2}{*}{PTDF Based Formulation}} & \multirow{2}{*}{} & \multicolumn{7}{c}{\multirow{2}{*}{Mixed Formulation}} \\
            \multicolumn{2}{c}{} & \multicolumn{7}{c}{} &  & \multicolumn{7}{c}{} \\ \cline{3-17} 
            \multicolumn{2}{c}{} & 20\% & 50\% & 70\% & 90\% & 95\% & 98\% & 99\% &  & 20\% & 50\% & 70\% & 90\% & 95\% & 98\% & 99\% \\ \hline
            \multirow{5}{*}{\begin{tabular}[c]{@{}c@{}}Delete\\ Random\\ Branches\end{tabular}} & case1354pegase & 6.40 & 5.48 & 3.64 & 1.33 & \textbf{0.60} & \textbf{0.46} & \textbf{0.38} &  & \textbf{1.48} & \textbf{1.16} & \textbf{1.13} & \textbf{0.80} & 0.70 & 0.63 & 0.63 \\
             & case1951rte & 9.16 & 6.11 & 3.44 & 1.20 & \textbf{0.61} & \textbf{0.44} & \textbf{0.40} &  & \textbf{1.76} & \textbf{1.59} & \textbf{1.26} & \textbf{1.05} & 0.86 & 0.84 & 0.80 \\
             & case3375wp & 25.92 & 24.44 & 17.89 & 3.21 & \textbf{1.44} & \textbf{0.70} & \textbf{0.52} &  & \textbf{5.54} & \textbf{6.80} & \textbf{7.27} & \textbf{3.04} & 3.73 & 2.48 & 2.50 \\
             & case6470rte & 398.86 & 183.44 & 67.05 & 9.02 & 4.53 & \textbf{2.70} & \textbf{2.39} &  & \textbf{7.16} & \textbf{6.71} & \textbf{5.05} & \textbf{4.20} & \textbf{4.17} & 4.11 & 3.96 \\
             & case6515rte & 256.62 & 165.53 & 72.80 & 11.91 & 4.52 & \textbf{2.79} & \textbf{2.02} &  & \textbf{7.34} & \textbf{5.56} & \textbf{4.97} & \textbf{3.78} & \textbf{3.50} & 3.25 & 3.33 \\ \hline
             \multirow{5}{*}{\begin{tabular}[c]{@{}c@{}}Delete\\ Uncongested\\ Branches\end{tabular}} & case1354pegase & 6.53 & 5.30 & 7.35 & 2.41 & 0.91 & \textbf{0.53} & \textbf{0.38} &  & \textbf{1.66} & \textbf{2.10} & \textbf{1.24} & \textbf{0.99} & 0.79 & 0.89 & 0.77 \\
             & case1951rte & 8.99 & 7.52 & 5.10 & 1.75 & \textbf{1.29} & \textbf{0.60} & \textbf{0.45} &  & \textbf{1.79} & \textbf{1.64} & \textbf{1.47} & \textbf{1.17} & 1.72 & 1.17 & 1.01 \\
             & case3375wp & 24.96 & 25.01 & 20.07 & 7.98 & \textbf{4.11} & \textbf{1.48} & \textbf{0.97} &  & \textbf{5.84} & \textbf{6.82} & \textbf{8.39} & \textbf{3.54} & 4.22 & 3.17 & 3.11 \\
             & case6470rte & 322.01 & 233.51 & 164.37 & 63.07 & 20.98 & 7.63 & \textbf{3.14} &  & \textbf{8.37} & \textbf{8.17} & \textbf{6.84} & \textbf{5.57} & \textbf{5.71} & \textbf{6.06} & 6.02 \\
             & case6515rte & 265.71 & 281.22 & 176.11 & 57.51 & 18.49 & 4.96 & \textbf{2.51} &  & \textbf{8.33} & \textbf{7.89} & \textbf{7.03} & \textbf{6.06} & \textbf{5.48} & \textbf{4.94} & 4.86 \\ \hline\hline
            \end{tabular}
      \begin{tablenotes}   
        \footnotesize              
        \item[*]  All the results are the average of five tests.
        \item[**]  Time is obtained using solver time (from calling the solver to getting the final result)
      \end{tablenotes} 
      \end{threeparttable}
      \label{tab:SCED_del_result}
\end{table*}
   
\begin{figure}[htbp!]
  \centering
  \includegraphics[width=0.45\textwidth]{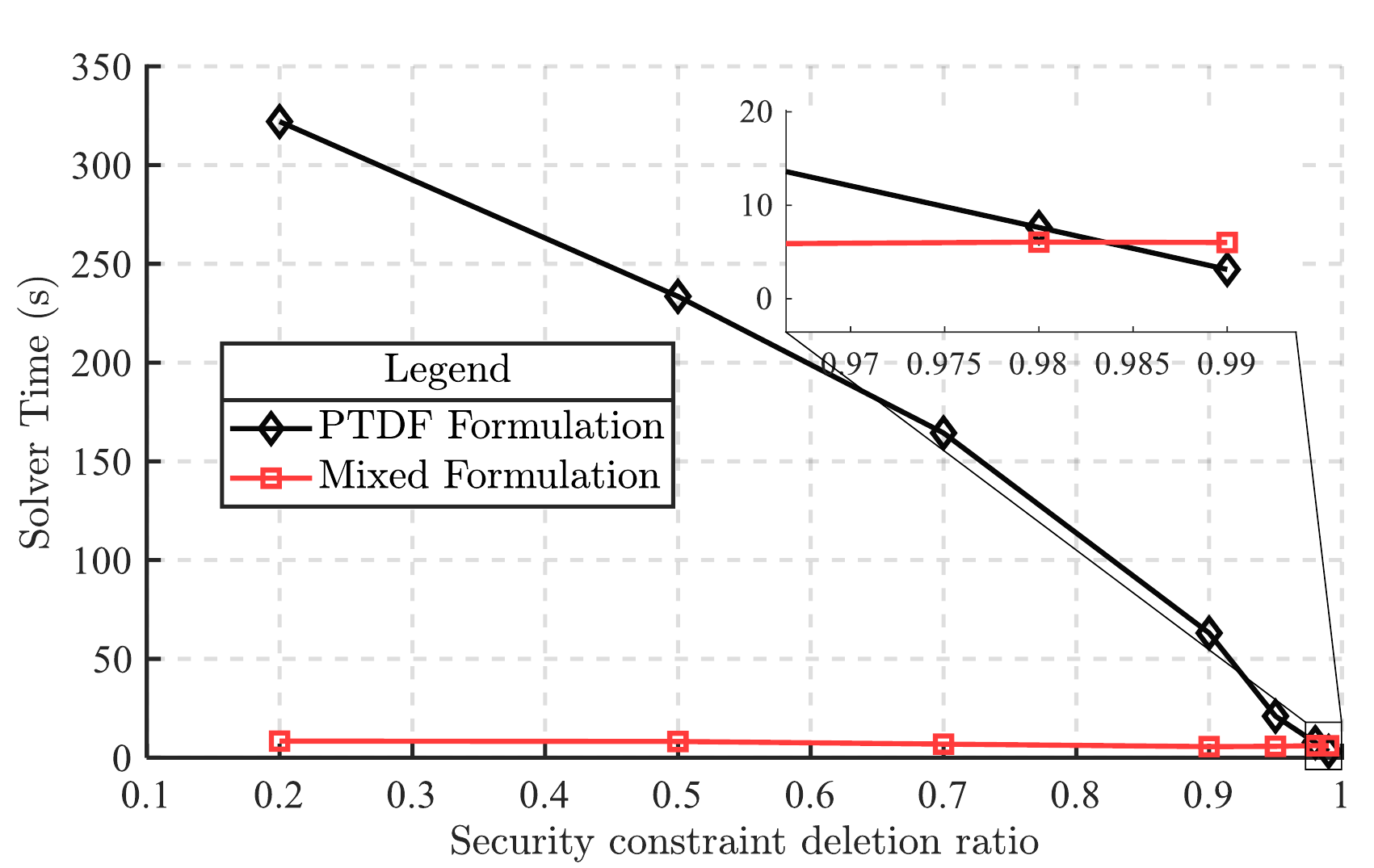}
  \caption{Solving time with varying the security constraints deletion ratio, using removing the most uncongested branches in case6470rte.}
  \label{fig:del}
\end{figure}

Under the two techniques, the solving time for both formulations decreases as the proportion of deleted branches increases. However, since the density of the PTDF formulation reduces with each deletion, the time required to solve the problem decreases more rapidly than that of the mixed formulation, whose density remains constant. Moreover, as the number of branch constraints reduces, the $P_{bus}=B_{bus}\theta$ constraint in the mixed formulation remains the same, resulting in the PTDF formulation requiring less time to solve than the mixed formulation once all branch constraints have been removed. As shown in Fig. \ref{fig:del}, using the approach of eliminating the most uncongested branch, the solving time of PTDF is less than that of the mixed formulation when the deletion ratio exceeds 98\%, while the mixed formulation is always superior before 98\%. According to the findings of reference \cite{zhai2010fast}, the fast identifying method can detect approximately 85\% of inactive security constraints.

\subsection{Discussion on constraints in mixed-integer problem}

In general, when solving SCUC, a solver typically employs a branch-and-bound or branch-and-cut framework, which sequentially adds linear constraints to the relaxation problem. IPM cannot use the last optimal point for a warm start, but the dual simplex technique can. Because the dual simplex approach is less sensitive to the density of the coefficients than IPM but more sensitive to the problem size, the PTDF formulation is always preferred over the mixed formulation for solving SCUC because it results in a smaller problem size.

The experimental results show that for a medium-sized case, if the default solution method, i.e., the dual simplex method, is used for the relaxation problem in the branch-and-bound framework, the PTDF formulation is significantly faster than the mixed formulation in terms of the node exploration speed and the convergence speed of the MIP gap. Moreover, if the IPM is used for the relaxation problem, it is difficult to find a feasible solution for several hours, and its solution speed is much slower than that of the dual simplex method.

\section{Conclusion}

In order to reduce the impact of the dense matrix generated by the PTDF formulation (pure generator power formulation) during each iteration of the primal-dual interior point method, we employ a matrix transformation technique to restore the sparsity of the linear system. Further, we find that this matrix transformation technique is equivalent to solving the problem rewritten with pure voltage phase angle formulation, and we obtain the mixed formulation (B-$\theta$ formulation) based on the pure voltage phase angle formulation. After extensive numerical experiments, the mixed formulation is similar to the pure voltage phase angle formulation in terms of solution speed, but the numerical stability is better. In the optimal power flow problems and security-constrained economic dispatch problems, the speed of the mixed formulation is faster than the PTDF formulation, and the speed advantage increases as the model size increases. In the SCUC problem, the PTDF formulation is recommended because reoptimization techniques are significantly more critical and the dual simplex method is typically employed rather than IPM in mixed integer problems. These conclusions hold true for linear programming formulations as well.


%





\ifCLASSOPTIONcaptionsoff
  \newpage
\fi



\bibliographystyle{IEEEtran}
\bibliography{IEEEabrv, ref.bib}
\end{document}